\documentclass{article}

\usepackage[top=1in,bottom=1in,left=1in,right=1in]{geometry}
\usepackage{amstext,amssymb,amsmath,amsthm,latexsym}
\usepackage{enumerate,bm}
\usepackage{mathrsfs,amsfonts}
\usepackage{paralist}
\usepackage[pdftex,colorlinks,citecolor=blue,linkcolor=blue]{hyperref}
\usepackage{algpseudocodex}
\usepackage[all]{xy}
\usepackage{longtable}

\theoremstyle{plain}
\newtheorem{theorem}{Theorem}[section]

\newtheorem{corollary}[theorem]{Corollary}
\newtheorem{lemma}[theorem]{Lemma}

\theoremstyle{definition}
\newtheorem{definition}[theorem]{Definition}

\newtheorem{example}[theorem]{Example}
\newtheorem{notation}[theorem]{Notation}

\theoremstyle{remark}
\newtheorem{remark}[theorem]{Remark}

\DeclareMathOperator*{\lcm}{lcm}
\DeclareMathOperator*{\rcm}{rcm}

\newcommand{\aar}{a}
\newcommand{\Bq}{\pb_\qr}
\newcommand{\bN}{{\BN^{n-1}}}
\newcommand{\bM}{\bm{[}\tx\bm{]}}
\newcommand{\BN}{\mathbb{N}}
\newcommand{\BQ}{\mathbb{Q}}
\newcommand{\BZ}{\mathbb{Z}}
\newcommand{\bg}{Buchberger}

\newcommand{\CS}{\mathcal{S}}
\newcommand{\cds}{\Omega}

\newcommand{\clm}{\tx^\gamma}
\newcommand{\cpel}{\pel^c}
\newcommand{\cps}{\Theta}
\newcommand{\dss}{\Lambda}

\newcommand{\el}{\chi}
\newcommand{\fd}{K}
\newcommand{\gb}{Gr\"obner}
\newcommand{\Gd}{\textnormal{\textsc{gcd}}}

\newcommand{\Iq}{I_\qr}
\newcommand{\Iu}{I\cap\kux}
\newcommand{\ibr}{b}
\newcommand{\icr}{c}
\newcommand{\icsp}{intermediate coefficient swell problem}
\newcommand{\idr}{d}

\newcommand{\ioq}{\iota_\qr}
\newcommand{\ipr}{p}
\newcommand{\iur}{\mu}
\newcommand{\kux}{K[x_1]}
\newcommand{\kxx}{\kux [\tx]}
\newcommand{\Lst}[3]{#1_{#2},\dotsc,#1_{#3}}
\newcommand{\lc}{\textrm{lc}}
\newcommand{\lcc}{\textnormal{\textsc{lc}}}
\newcommand{\Lcm}{\textnormal{\textsc{lcm}}}
\newcommand{\lex}{\textnormal{\textsc{lex}}}
\newcommand{\lm}{\textrm{lm}}
\newcommand{\lmc}{\textnormal{\textsc{lm}}}
\newcommand{\lmr}{m}
\newcommand{\lnr}{n}
\newcommand{\lt}{\textrm{lt}}
\newcommand{\ltc}{\textnormal{\textsc{lt}}}
\newcommand{\Mo}{M\"oller}
\newcommand{\mar}{\lambda}
\newcommand{\mas}{\Lambda}
\newcommand{\mel}{\chi_\qr}
\newcommand{\mlt}{\mathrm{mult}}
\newcommand{\mop}{modular proper}
\newcommand{\mq}{\sigma_\qr}

\newcommand{\nonk}{\kux\setminus K}
\newcommand{\Pid}{principal ideal domain}
\newcommand{\Pir}{principal ideal ring}

\newcommand{\pb}{B}
\newcommand{\pel}{\chi_0}
\newcommand{\pid}[1]{(#1)}
\newcommand{\prb}{proper basis}
\newcommand{\qr}{q}

\newcommand{\RD}{\kxx\setminus\kux}
\newcommand{\Rq}{R_\qr}
\newcommand{\Rqx}{\Rq [\tx]}
\newcommand{\Rqd}{\Rqx\setminus\Rq}

\newcommand{\rat}{\Rq^\ast\setminus\Rq^\times}
\newcommand{\rnd}{(\Rqx)^\ast\setminus\Rq^\times}

\newcommand{\rr}{r}
\newcommand{\Sp}{$S$-polynomial}
\newcommand{\supp}[1]{\mathrm{supp}(#1)}

\newcommand{\tM}{\bm{[}\tx\bm{]}}
\newcommand{\tas}{F}
\newcommand{\tbs}{G}
\newcommand{\tf}{f}
\newcommand{\tg}{g}
\newcommand{\tx}{\tilde{\X}}
\newcommand{\X}{\bm{x}}
\newcommand{\zpi}{zero-dimensional polynomial ideal}

\numberwithin{equation}{section}

\begin{document}

\title{The Proper Basis for Polynomial Ideals
\footnotetext{\textit{Email:} \texttt{smmath@foxmail.com; masm@buaa.edu.cn}}
\footnotetext{\textit{Address:} School of Mathematical Sciences, Beihang University, Beijing 100191, China.}
\footnotetext{2020 \textit{Mathematics Subject
Classification.} Primary 13P10; 13B25.}
\footnotetext{\textit{Key words and phrases:} polynomial ideal, ideal basis, \prb, \gb\ basis, proper division, \icsp\ (ICSP).}}
\author{Sheng-Ming Ma}
\date{}

\maketitle

\begin{abstract}
We define a new type of ideal basis called the \prb\ that improves both \gb\ basis and \bg's algorithm.
Let $x_1$ be the least variable of a monomial ordering in a polynomial ring $K[\Lst x1n]$ over a field $K$.
The \gb\ basis of a \zpi\ contains a univariate polynomial in $x_1$.
The \prb\ is defined and computed in the variables $\tx:=(\Lst x2n)$ with $x_1$ serving as a parameter in the algebra $\kxx$.
Its algorithm is more efficient than not only \bg's algorithm whose elimination of $\tx$ unnecessarily involves the least variable $x_1$ but also \Mo's algorithm due to its polynomial division mechanism.
This is corroborated by a series of benchmark testings herein.
The \prb\ is in a modular form and neater than \gb\ basis and hence reduces its coefficient swell problem.
It is expected that all the state of the art algorithms improving \bg's algorithm over the last decades can be further improved if we apply them to the proper basis.
\end{abstract}

\section{Introduction}
After \bg\ initiated his celebrated algorithm in his remarkable PhD thesis \cite{Buc65}, the theory of \gb\ basis had been established as a standard tool in computer algebra, yielding algorithmic solutions to many significant problems in mathematics, science and engineering \cite{BW98}.
As a result, there have been many excellent textbooks on the subject such as \cite{BW93,Mis93,AL94,KR00,GP08,EH12,GG13,CLO15}.

The computation of \gb\ basis is often plagued with a high complexity.
A typical phenomenon is the \icsp\ (ICSP) over the rational number field $\BQ$ and especially with respect to the \lex\ ordering.
This stimulates decades of ardent endeavors to improve the efficiency of \bg's algorithm.
The methods of normal selection strategies and signatures have been quite successful in this respect \cite{Buc85,GMN91,BW93,Fau02,GVW16,EF17,FV20}.
The modular and $p$-adic methods based on the ``lucky primes" and Hensel lifting are adopted to control the rampant growth of the intermediate coefficients \cite{Ebe83,Win88,ST89,Tra89,Pau92,Gra93,Arn03}.
There are also the conversion methods among \gb\ bases such as the FGLM algorithm \cite{FGLM93} and \gb\ Walk \cite{CKM97}, a detailed description of which can be found in \cite{CLO05,Stu95}.

The \gb\ basis over a field had also been generalized to over rings.
In particular, the \gb\ basis over a \Pir\ (PIR) was developed by \Mo\ in \cite{Mol88} and elucidated in \cite[Chapter 4]{AL94}.
Nonetheless these generalizations have never been applied to polynomial ideals over a field.

The method of characteristic set \cite{Mis93,Wu01} based on the pseudo-division is more efficient than \gb\ basis.
However the pseudo-division usually loses too much algebraic information of the original ideal.
After the computation we only have information on the zero locus or radical ideal that is insufficient to solve algebraic problems.
This is also the deficiency associated with a few other methods such as the rational univariate representation \cite{Rou99}.

In this article we define a new type of ideal basis called \prb\ that reduces the computational complexity of \gb\ basis while retains the algebraic information of the original ideal.
Let $I\subset K[\X]$ be a \zpi\ over a field $K$ with $\X:=(\Lst x1n)$.
Suppose that $\Iu$ is a principal ideal $\pid\el$ that is generated by the eliminant $\el\in\nonk$ of $I$.
We obtain $\el$ after eliminating all the other variables $\tx:=(\Lst x2n)$ with respect to, e.g., the \lex\ ordering $x_1\prec\cdots\prec x_n$.
Nonetheless it is literally excessive computations for \bg's elimination process to involve the variable $x_1$.
Hence it is natural for us to treat $K[\X]$ as the algebra $\kxx$ such that the variable $x_1$ serves as a parameter when we eliminate the variables $\tx$.
In this way we define a new type of ideal basis for $I$ called the \prb.
Moreover, the \prb\ is based on a new type of polynomial division called the proper division, which improves the division mechanism in \Mo's algorithm over a PIR for \gb\ basis.
Our benchmark testings corroborate that the \prb\ algorithm is distinctively more efficient than both \bg's classical algorithm over $\fd$ and \Mo's one over $\kux$ for \gb\ bases.

In Definition \ref{Def:TermReduction} and Theorem \ref{Thm:PseudoReduction} we define the proper division using least multipliers in $\kux$.
The purpose of Lemma \ref{Lemma:RelativePrimePairs} and Lemma \ref{Lemma:TriangleIdentity} is to trim down the number of \Sp s for computational efficiency.
The algorithm in Lemma \ref{Algo:PseudoEliminant} parallels that of \bg\ for \gb\ basis except that it uses the least multiplier in $\kux$.
Based on Lemma \ref{Lemma:PseudoReps}, we derive in Corollary \ref{Cor:LeadTermDef} \eqref{item:ComptblProperBas} that the prebasis obtained in Lemma \ref{Algo:PseudoEliminant} satisfies Definition \ref{Def:NewBases} for the \prb.

The content of Section \ref{Sec:ModularAlgm} parallels Section \ref{Sec:DivisionAlgm} except that we contrive a modular algorithm over a PIR $\kux/\pid\qr$ that might contain zero divisors.

In Section \ref{Sec:Timings} we conduct benchmark testings on the timings of the respective algorithms.
It is clear that the \prb\ algorithm has a distinctive advantage over both \bg's and \Mo's classical algorithms for \gb\ bases in the textbooks.

Throughout the article all our discussions are with respect to the \lex\ ordering $x_1\prec\cdots\prec x_n$ since it is typical to have the highest level of computational complexity compared with other monomial orderings.
We use $K$ to denote a field that is not necessarily algebraically closed unless specified.
As usual, let us denote the sets of rational, integral and natural numbers by $\BQ$, $\BZ$ and $\BN$ respectively.

\section{The Definition for Proper Basis}\label{Sec:InductiveGroebner}

Let $K$ be a field and $\X$ denote the variables $(\Lst x1n)$.
We consider two types of algebras that determine the notations in this article.
The first one is $\kxx$ over $\kux$ with the variables $\tx$ denoting $(\Lst x2n)$;
The second one is $\Rqx$ over a PIR $\Rq$ that is isomorphic to $\kux/\pid\qr$.
Here $\pid\qr$ denotes a nontrivial principal ideal in $\kux$ that is generated by $\qr$.

In order to clarify the notations, in what follows let us use $R$ to denote $\kux$ or $\Rq$.
In this way the algebra $R[\tx]$ over $R$ denotes the above algebra $\kxx$ or $\Rqx$.
We also adopt the following notations for a ring $R$:
we denote $R^\ast:=R\setminus\{0\}$ as the set of nonzero elements in $R$, and $R^\times$ as the set of units in $R^\ast$.

With $\alpha=(\Lst\alpha 2n)\in\bN$, we denote a \emph{monomial} $x_2^{\alpha_2}\cdots x_n^{\alpha_n}$ by $\tx^\alpha$ and a \emph{term} by $c\tx^\alpha$ with the \emph{coefficient} $c\in R^\ast$.
Let us denote the set of monomials in $\tx$ by $\bM:=\{\tx^\alpha\colon\alpha\in\bN\}$.
The notation $\langle A\rangle$ denotes an ideal generated by a nonempty subset $A\subset R[\tx]$.

\begin{notation}[$\supp f$; $\ltc(f)$, $\lmc (f)$, $\lcc (f)$; $\lmc (\pb)$, $\langle\lmc (\pb)\rangle$, $\ltc (\pb)$, $\langle\ltc (\pb)\rangle$]\label{Notation:LeadingEntities}
\hfill

Let $R$ denote a PIR such as $\kux$ or $\Rq\cong\kux/\pid\qr$ as above.
Suppose that $f=\sum_\alpha c_\alpha\tx^\alpha$ is a polynomial in $R[\tx]$.
We denote the \emph{support} of $f$ by $\supp f:=\{\tx^\alpha\in\bM\colon c_\alpha\ne 0\}$.
In particular, we define $\supp f:=\{1\}$ if $f\in R^\ast$ and $\supp f:=\emptyset$ if $f=0$.

Let $\succ$ be a monomial ordering on the monomial set $\bM$.
The \emph{leading term} of $f$ is a term $c_\beta\tx^\beta$ that satisfies $\tx^\beta:=\max_{\succ}\{\tx^\alpha\in\supp f\}$ and is denoted by $\ltc (f):=c_\beta\tx^\beta$.
The \emph{leading monomial} of $f$ is the monomial $\tx^\beta$ of $\ltc (f)=c_\beta\tx^\beta$ and is denoted by $\lmc (f):=\tx^\beta$.
The \emph{leading coefficient} of $f$ is the coefficient $c_\beta$ of the leading term $\ltc (f)=c_\beta\tx^\beta$ and is denoted by $\lcc (f):=c_\beta\in R^\ast$.

Let $\pb=\{\ibr_j\colon 1\le j\le s\}$ be a polynomial set in $R[\tx]\setminus\{0\}$.
We denote the leading monomial set $\{\lmc (\ibr_j)\colon 1\le j\le s\}$ by $\lmc (\pb)$, and the ideal generated by $\lmc (\pb)$ in $R[\tx]$ by $\langle\lmc (\pb)\rangle$.
Similarly we can define $\ltc (\pb)$ and $\langle\ltc (\pb)\rangle$.
\qed\end{notation}

For $\qr\in\nonk$, consider the set $\Rq:=\{r\in\kux\colon\deg (r)<\deg (\qr)\}$ with $\deg (r)=0$ for every $r\in K$ including $r=0$.
We define binary operations on $\Rq$ such that it is a PIR isomorphic to $\kux/\pid\qr$.
For every $f\in\kux$, there exist a quotient $h\in\kux$ and unique remainder $r\in\Rq$ such that $f=h\qr+r$.
We define an epimorphism as follows.
\begin{equation}\label{ProjectionPQR}
\mq\colon\kux\rightarrow\Rq\colon\quad\mq (f):=r.
\end{equation}

Since $\Rq$ is a subset of $\kux$, for every $r\in\Rq$, we define an injection as follows.
\begin{equation}\label{PQREmbedding}
\ioq\colon\Rq\hookrightarrow\kux\colon\quad\ioq (r):=r.
\end{equation}

If we extend the ring epimorphism $\mq$ in \eqref{ProjectionPQR} such that it is the identity map on the variables $\tx$, then $\mq$ induces an epimorphism from $\kxx$ to $\Rqx$ which we still denote by $\mq$ as follows.
\begin{equation}\label{ExtendedProjectionPQR}
\mq\colon\kxx\rightarrow\Rqx\colon\quad\mq\Bigl(\sum_{j=1}^s c_j\tx^{\alpha_j}\Bigr):=\sum_{j=1}^s\mq (c_j)\tx^{\alpha_j}.
\end{equation}

Similarly the injection $\ioq$ in \eqref{PQREmbedding} can be extended to an injection from $\Rqx$ into $\kxx$ in the way that it is an identity map on the variables $\tx$ as follows.
\begin{equation}\label{ExtendedEmbedding}
\ioq\colon\Rqx\hookrightarrow\kxx\colon\quad\ioq\Bigl(\sum_{j=1}^s c_j\tx^{\alpha_j}\Bigr):=\sum_{j=1}^s\ioq(c_j)\tx^{\alpha_j}.
\end{equation}

\begin{definition}[Eliminant $\el$; multiplicity $\mlt_\ipr (f)$]\label{Def:Eliminant}
\hfill

For a \zpi\ $I\subset\kxx$, let us denote the monic generator of the principal ideal $\Iu$ by $\el$, i.e., $\Iu=\pid\el$.
We call $\el$ the \emph{eliminant} of $I$ henceforth.

For an irreducible factor $\ipr$ of a univariate polynomial $f\in\nonk$,  we use $\mlt_\ipr (f)$ to denote the multiplicity of $\ipr$ in $f$.
That is, $\mlt_\ipr (f):=\max\{i\in\BN\colon f\in\pid{\ipr^i}\}$.
\qed\end{definition}

Suppose that the above eliminant $\el$ has the following factorization whose pairwise relatively prime factors are in $\cds\subset\nonk$:
\begin{equation}\label{EliFactors}
\el=\prod_{\qr\in\cds}\qr.
\end{equation}
It is evident that the above factorization corresponds to a decomposition of $I$ as follows.
\begin{equation}\label{IdealDecomposition}
I=\bigcap_{\qr\in\cds}(I+\langle\qr\rangle).
\end{equation}

\begin{definition}[Proper basis]\label{Def:NewBases}
\hfill

Let $I\subset\kxx$ be a \zpi\ whose eliminant $\el$ has a pairwise coprime factorization as in \eqref{EliFactors}.
For each $\qr\in\cds$, let $\Iq:=\mq (I)\subset\Rqx$ denote its image under $\mq$ in \eqref{ExtendedProjectionPQR}.
Also suppose that there exists a subset $\Bq\subset\Iq$ that satisfies:
\begin{equation}\label{ProperLeadGen}
\langle\ltc (\Iq)\rangle=\langle\ltc (\Bq)\rangle.
\end{equation}
Then we define the following set:
\begin{equation}\label{NewBases}
\bigcup_{\qr\in\cds}\bigl(\Bq\cup\{\qr\}\bigr)
\end{equation}
as the \emph{\prb} of the \zpi\ $I$.
\end{definition}

\section{The First Algorithmic Step in $\kxx$}\label{Sec:DivisionAlgm}

Our algorithm for the \prb\ is divided into two steps.
The first step is in $\kxx$ over $\kux$ whereas the second one is in $\Rqx$ over $\Rq$ as in Section \ref{Sec:ModularAlgm}.

\begin{notation}[$\gcd (a,b)$, $\lcm (a,b)$]
\hfill

For a univariate polynomial $f=\sum_{k=0}^{d}c_kx_1^k\in\kux$ with $c_d\in\fd^\ast$ being nonzero, recall that we use $\lt (f)$, $\lm (f)$ and $\lc (f)$ to denote the leading term $c_dx_1^d$, leading monomial $x_1^d$ and leading coefficient $c_d$ of $f$ over the field $\fd$ respectively.
In what follows we use $\gcd (a,b)$ and $\lcm (a,b)$ to denote the greatest common divisor (\Gd) and least common multiple (\Lcm) of $a,b\in (\kux)^\ast$ respectively.
In particular, we always choose the monic polynomials for \Gd s such that $\lc (\gcd (a,b))=1$.
However we always choose the \Lcm s such that $\lc (\lcm (a,b))=\lc (a)\cdot\lc (b)$.
This is a different choice from that of \cite[\S 3.4]{GG13} so that the identity $a/\gcd (a,b)=\lcm (a,b)/b$ holds, which is convenient for our later discussions.
\end{notation}

\begin{definition}[Proper term reduction in the algebra $\kxx$]\label{Def:TermReduction}
\hfill

Let $\succ$ be a monomial ordering on $\bM$.
For $f\in\RD$ and $g\in\kxx\setminus\{0\}$, suppose that $f$ has a term $c_\alpha\X^\alpha$ satisfying $\X^\alpha\in\supp f\cap\langle\lmc (g)\rangle$.
We define a \emph{proper term reduction} of $c_\alpha\X^\alpha$ by $g$ as follows.
\begin{equation}\label{TermReduction}
h=\iur f-\frac{\lmr\X^\alpha}{\ltc (g)}g
\end{equation}
with the \emph{least multiplier} $\iur:=\lmr/c_\alpha\in R^\ast$ and multiplier $\lmr:=\lcm (c_\alpha,\lcc (g))$.
\qed\end{definition}

\begin{definition}[Properly reduced polynomial in $\kxx$]
\hfill

A polynomial $\rr\in\kxx$ is \emph{properly reduced} with respect to a polynomial set $\pb=\{\ibr_j\colon 1\le j\le s\}\subset\RD$ if $\supp\rr\cap\langle\lmc (\pb)\rangle=\emptyset$ holds.
In particular, this includes the special case when $\rr=0$ and hence $\supp\rr=\emptyset$.
Otherwise we say that $\rr$ is \emph{properly reducible} with respect to $\pb$.
\qed\end{definition}

\begin{theorem}[Proper division or reduction in $\kxx$]\label{Thm:PseudoReduction}
\hfill

Suppose that $\pb=\{\ibr_j\colon 1\le j\le s\}\subset\RD$ is a finite polynomial set.
For every $f\in\kxx$, there exist a multiplier $\mar\in (\kux)^\ast$, a remainder $\rr\in\kxx$ and quotients $\qr_j\in\kxx$ for $1\le j\le s$ such that
\begin{equation}\label{PseudoDivisionExpression}
\mar f=\sum_{j=1}^s\qr_j\ibr_j+\rr,
\end{equation}
where $\rr$ is properly reduced with respect to $\pb$, and the multiplier $\mar$ is literally a product of the least multipliers as in \eqref{TermReduction}.
Moreover, the polynomials in \eqref{PseudoDivisionExpression} satisfy the following condition:
\begin{equation}\label{DivisionCond}
\lmc (f)=\max\bigl\{\max_{1\le j\le s}\{\lmc (\qr_j\ibr_j)\},\lmc (\rr)\bigr\}.
\end{equation}
\end{theorem}
\begin{proof}
If $f$ is properly reducible with respect to $\pb$, we define $\X^\alpha:=\max_{\succ}\{\supp f\cap\langle\lmc (\pb)\rangle\}$ and make a proper term reduction of $c_\alpha\X^\alpha$ by some $\ibr_j\in\pb$ as in \eqref{TermReduction}.
This repeats and then terminates in finite steps since the monomial ordering is a well-ordering.
The equality in \eqref{DivisionCond} follows from the fact that it holds for the proper term reduction in \eqref{TermReduction}.
\end{proof}

\begin{remark}\label{Rmk:LinearEqs}
There is a stark difference between the proper reduction herein and \Mo's reduction over a PIR in \cite{Mol88}.
With the notations as in Theorem \ref{Thm:PseudoReduction}, \Mo's reduction requires that the linear equation
\begin{equation}\label{LinearCoeffEq}
\lcc (f)=\sum_{j=1}^s\icr_j\cdot\lcc (\ibr_j)
\end{equation}
be solvable for the $\icr_j$'s in $\kux$.
Please refer to \cite[P349, (1)]{Mol88} or \cite[P204, (4.1.1)]{AL94}.
This is the major difference between the proper reduction and \Mo's reduction.
\qed\end{remark}

\begin{definition}[\Sp\ in $\kxx$]\label{Def:SPolynomial}
\hfill

Suppose that $f,g\in\RD$.
The \emph{\Sp} of $f$ and $g$ is defined as:
\begin{equation}\label{SPolynomialDef}
S(f,g):=\frac{\lmr\clm}{\ltc (f)}f-\frac{\lmr\clm}{\ltc (g)}g
\end{equation}
with $\lmr:=\lcm (\lcc (f),\lcc (g))\in (\kux)^\ast$ and $\clm:=\lcm (\lmc (f),\lmc (g))\in\tM$.

In particular, when $g\in (\kux)^\ast$ and $f\in\RD$, we take $\lmc (g)=1$ and $\lmr=\lcm (\lcc (f),g)$ and the special \Sp\ satisfies $\idr S(f,g)=f_1g$ with $\idr:=\gcd (\lcc (f),g)\in\kux$ and $f_1:=f-\ltc (f)$.
\qed\end{definition}

It is easy to verify that the \Sp\ satisfies $\lmc (S(f,g))\prec\clm=\lcm (\lmc (f),\lmc (g))$.

Lemma \ref{Lemma:RelativePrimePairs} and Lemma \ref{Lemma:TriangleIdentity} are generalizations of \bg's first and second criteria respectively.
Please refer to \cite[P222, \S 5.5]{BW93} or \cite[P124, \S 3.3]{AL94} for \bg's two criteria.

\begin{lemma}\label{Lemma:RelativePrimePairs}
For $f,g\in\RD$, suppose that $\lmc (f)$ and $\lmc (g)$ are relatively prime.
Let us denote $\idr:=\gcd (\lcc (f),\lcc (g))$.
Then the \Sp\ in \eqref{SPolynomialDef} satisfies $\idr S(f,g)=f_1g-g_1f$ with $f_1:=f-\ltc (f)$ and $g_1:=g-\ltc (g)$.
Moreover, we have $\lmc (S(f,g))=\max\{\lmc (f_1g),\lmc (g_1f)\}$.
\end{lemma}
\begin{proof}
With $\lmr$ denoting $\lcm (\lcc (f),\lcc (g))$, the conclusion follows from the identities $\lmr/\lcc (f)=\lcc (g)/\idr$ and $\lmr/\lcc (g)=\lcc (f)/\idr$ .
Moreover, the identity for $\lmc (S(f,g))$ follows from $\idr S(f,g)=f_1\cdot\ltc (g)-g_1\cdot\ltc (f)$ that is easy to prove.
\end{proof}

Let us denote $\lcm (\ltc (f),\ltc (g)):=\lcm (\lcc (f),\lcc (g))\cdot\lcm (\lmc (f),\lmc (g))$.
The following lemma is evident.

\begin{lemma}\label{Lemma:TriangleIdentity}
For $f,g,h\in\RD$, if $\lcm (\lmc (f),\lmc (g))$ is divisible by $\lmc (h)$, then we have the following triangular identity among the \Sp s:
\begin{equation}\label{TriangleIdentity}
\mar S(f,g)=\frac{\mar\cdot\lcm (\ltc (f),\ltc (g))}{\lcm (\ltc (f),\ltc (h))}S(f,h)-\frac{\mar\cdot\lcm (\ltc (f),\ltc (g))}{\lcm (\ltc (g),\ltc (h))}S(g,h),
\end{equation}
where the multiplier is defined as $\mar:=\lcc (h)/\idr$ with $\idr:=\gcd (\lcm (\lcc (f),\lcc (g)),\lcc (h))\in\kux$.
\end{lemma}

\begin{lemma}[The algorithm computing a temporary prebasis for the \prb]\label{Algo:PseudoEliminant}
\hfill

Let $\tas\subset\RD$ be a finite polynomial set generating a \zpi.
Then the algorithm that consists of the following sequence of instructions terminates in finite steps.
\begin{enumerate}
\renewcommand{\labelenumi}{\textbf{\upshape{A\theenumi}}}

  \item Start with a prebasis set $\tbs:=\tas$, multiplier set $\mas:=\emptyset$ in $\kux$, temporary set $\CS:=\emptyset$ of \Sp s, and temporary pre-eliminant $\pel:=0$.

  \item For each pair $\tf,\tg\in\tbs$ with $\tf\ne\tg$, execute the instructions in \textbf{\upshape{A\ref{item:SPolyCompt}}}.

   \item\label{item:SPolyCompt} If $\lmc (\tf)$ and $\lmc (\tg)$ are relatively prime, add $\idr:=\gcd (\lcc (f),\lcc (g))$ into the multiplier set $\mas$ if $\idr\in\nonk$;

   Else if there exists an $h\in\tbs\setminus\{\tf,\tg\}$ such that $\lcm (\lmc (f),\lmc (g))\in\langle\lmc (h)\rangle$, add $\mar$ as defined in \eqref{TriangleIdentity} into the multiplier set $\mas$ if $\mar\in\nonk$;

   Else, compute the \Sp\ $S(\tf,\tg)$ as in \eqref{SPolynomialDef} and add it into the set $\CS$.

  \item\label{item:SPolyReduct} For every \Sp\ $S\in\CS$, make a proper reduction of $S$ by the prebasis set $\tbs$ to obtain a remainder that is denoted as $r$ like in Theorem \ref{Thm:PseudoReduction}.
   If the multiplier $\mar$ satisfies $\mar\in\nonk$ in the proper reduction \eqref{PseudoDivisionExpression}, add $\mar$ into the multiplier set $\mas$.
   Delete $S$ from the set $\CS$.

   If $\rr\in\nonk$, redefine $\pel:=\gcd (\rr,\pel)$;

   Else if $\rr\in\RD$, add $\rr$ into $\tbs$.
   For every $\tf\in\tbs\setminus\{\rr\}$, execute the instructions in
   \textbf{\upshape{A\ref{item:SPolyCompt}}} to compute the \Sp\ $S(\tf,\rr)$;

   Else if $\rr\in K^\ast$ is a unit, halt the algorithm and output $\tbs=\{1\}$.

 \item Repeat the instructions in \textbf{\upshape{A\ref{item:SPolyReduct}}} until $\CS=\emptyset$.

 \item\label{item:SpecialReduct} For every $\tf\in\tbs$, if $\idr:=\gcd (\lcc (f),\pel)\in\nonk$, add $\idr$ into the multiplier set $\mas$.

 \item Return the pre-eliminant $\pel\in (\kux)^\ast$, prebasis $\tbs\subset\RD$, multiplier set $\mas\subset\nonk$.
\end{enumerate}
\end{lemma}
\begin{proof}
The termination of the algorithm follows from the fact that the algebra $\kxx$ is Noetherian.
\end{proof}

Please note the step \textbf{\upshape{A\ref{item:SpecialReduct}}} is for the reduction of the special \Sp\ $S(\tf,\pel)$ as in Definition \ref{Def:SPolynomial}.

\begin{notation}[Compatible part $\cpel$ of the pre-eliminant $\pel$]\label{Nott:ComptblPart}
\hfill

Let $\pel$ and $\mas$ be the pre-eliminant and multiplier set that are obtained in the algorithm of Lemma \ref{Algo:PseudoEliminant}.
Let us denote $\cps:=\{\ipr\in\nonk\colon\ipr$ is an irreducible factor of $\pel$ and relatively prime to every multiplier in $\mas\}$.
If $\cps\ne\emptyset$, we define the \emph{compatible part} of $\pel$ as $\cpel:=\prod\limits_{\ipr\in\cps}\ipr^{\mlt_\ipr (\pel)}$;
Otherwise $\cpel:=1$.
\end{notation}

\begin{lemma}\label{Lemma:SyzygyTransform}
Let $\tas=\{f_j\colon 1\le j\le s\}\subset\RD$ be a finite polynomial set.
Suppose that each $f_j$ has the same leading monomial $\lmc (f_j)=\tx^\alpha\in\tM$ for $1\le j\le s$.
\begin{enumerate}[(1)]
\item\label{item:SPolynomialExpansion} If $f=\sum_{j=1}^s f_j$ satisfies $\lmc (f)\prec\tx^\alpha$, then
there exist multipliers $\ibr,\ibr_j\in (\kux)^\ast$ for $1\le j<s$ such that
\begin{equation}\label{SPolynomialExpansion}
\ibr f=\sum_{1\le j<s}\ibr_jS(f_j,f_s).
\end{equation}

\item\label{item:Nondivisibility} For each irreducible polynomial $\ipr\in\nonk$, we can always relabel the subscripts of the polynomial set $\tas=\{f_j\colon 1\le j\le s\}$ such that the multiplier $\ibr\in (\kux)^\ast$ in \eqref{SPolynomialExpansion} is not divisible by $\ipr$.
\end{enumerate}
\end{lemma}
\begin{proof}
\begin{inparaenum}[(1)]
\item Let us denote $l_j:=\lcc (f_j)$ for $1\le j\le s$ and $\lmr_j:=\lcm (l_j,l_s)$ for $1\le j<s$.
Then the condition $\lmc (f)\prec\tx^\alpha$ implies that $\sum_{j=1}^s l_j=0$.
It is easy to check that the identity \eqref{SPolynomialExpansion} holds for the definition of multipliers as $\ibr:=\lcm\{\lmr_j/l_j\colon 1\le j<s\}$ and $\ibr_j:=\ibr l_j/\lmr_j$ for $1\le j<s$.

\item Let us relabel the subscripts of $f_j$ for $1\le j\le s$ such that $\mlt_\ipr (l_s)=\min_{1\le j\le s}\{\mlt_\ipr (l_j)\}$.
Then $\mlt_\ipr (\lmr_j/l_j)=\mlt_\ipr (l_s/\gcd (l_j,l_s))=0$ for $1\le j<s$.
Thus the multiplier $\ibr$ is not divisible by $\ipr$.
\end{inparaenum}
\end{proof}

\begin{lemma}\label{Lemma:PseudoReps}
For a finite polynomial set $\tas\subset\RD$ generating a \zpi, let $\tbs=\{g_k\colon 1\le k\le s\}$ and $\pel$ be the prebasis and pre-eliminant that are obtained in the algorithm of Lemma \ref{Algo:PseudoEliminant}.
For every $f\in\langle\tas\rangle$, there exist $\{\qr_k\colon 0\le k\le s\}\subset\kxx$ and a multiplier $\mar\in (\kux)^\ast$ being relatively prime to the compatible part $\cpel$ of $\pel$ such that:
\begin{equation}\label{PseudoReps}
\mar f=\sum_{k=1}^s\qr_k\tg_k+\qr_0\pel.
\end{equation}
Moreover, the polynomials in \eqref{PseudoReps} satisfy the following condition:
\begin{equation}\label{PseudoRepCond}
\lmc (f)=\max\bigl\{\max_{1\le k\le s}\{\lmc (\qr_k\tg_k)\},\lmc (\qr_0)\bigr\}.
\end{equation}
In particular, the eliminant $\el\in\langle\tas\rangle$ also satisfies \eqref{PseudoReps} with $\qr_0\in(\kux)^\ast$ and $\gcd (\mar,\cpel)=1$ as follows:
\begin{equation}\label{FinalReduct}
\mar\el=\qr_0\pel.
\end{equation}
\end{lemma}
\begin{proof}
Let us fix an irreducible factor $\ipr$ of the compatible part $\cpel$, i.e., $\ipr\in\cps$ as in Notation \ref{Nott:ComptblPart}.

For $\tas=\{f_k\colon 1\le k\le t\}$, suppose that $f$ can be written as $f=\sum_{k=1}^th_kf_k$ with $h_k\in\kxx$ for $1\le k\le t$ and $\tx^\beta:=\max_{1\le k\le t}\{\lmc (h_kf_k)\}$.
Assume that $\lmc (f)\prec\tx^\beta$ since the conclusion already holds otherwise.
Let us denote $\ltc (h_k):=c_k\tx^{\alpha_k}$ with $c_k\in (\kux)^\ast$ for $1\le k\le t$.
Without loss of generality, suppose that $\lmc (\tx^{\alpha_k}f_k)=\tx^\beta$ holds for $1\le k\le t$.
And we have:
\begin{equation}\label{SepLeadTerm}
f=\sum_{k=1}^tc_k\tx^{\alpha_k}f_k+\sum_{k=1}^t(h_k-\ltc (h_k))f_k.
\end{equation}

As per Lemma \ref{Lemma:SyzygyTransform} \eqref{item:SPolynomialExpansion}, there exist multipliers $\ibr,\ibr_k\in (\kux)^\ast$ for $1\le k<t$ that satisfy the identity:
\begin{equation}\label{InvokeSPoly}
\ibr\sum_{k=1}^tc_k\tx^{\alpha_k}f_k=\sum_{1\le k<t}\ibr_kS(c_k\tx^{\alpha_k}f_k,c_t\tx^{\alpha_t}f_t)
=\sum_{1\le k<t}\ibr_k\lmr_k\tx^{\beta-\gamma_k}S(f_k,f_t),
\end{equation}
where $\tx^{\gamma_k}:=\lcm (\lmc (f_k),\lmc (f_t))$ and $\lmr_k:=\lcm (c_k\cdot\lcc (f_k),c_t\cdot\lcc (f_t))/\lcm (\lcc (f_k),\lcc (f_t))$ for $1\le k<t$.
Moreover, by Lemma \ref{Lemma:SyzygyTransform} \eqref{item:Nondivisibility}, we can relabel the subscripts of the elements in $\tas$ such that $\mlt_\ipr (\ibr)=0$.

Each \Sp\ $S(f_k,f_t)$ in \eqref{InvokeSPoly} has a proper reduction by the prebasis set $\tbs$ as in the steps \textbf{\upshape{A\ref{item:SPolyReduct}}} of the algorithm in Lemma \ref{Algo:PseudoEliminant}.
Further, we have $\lmc (S(\tf_k,\tf_t))\prec\tx^{\gamma_k}$.
These proper reductions together with the second summation in \eqref{SepLeadTerm} yield a new representation $\iur\ibr f=\sum_{k=0}^s\aar_kg_k$ with $g_0:=\pel$ such that $\max_{0\le k\le s}\{\lmc (\aar_kg_k)\}\prec\tx^\beta$.
Here the quotient $\aar_k$ is in $\kxx$ for $0\le k\le s$.
The multiplier $\iur\in (\kux)^\ast$ satisfies $\gcd (\iur,\cpel)=1$.
It is evident that $\mlt_\ipr (\iur\ibr)=0$.

We repeat the above procedure on the new representation $\iur\ibr f=\sum_{k=0}^s\aar_kg_k$.
The repetition halts in finite steps since the monomial ordering on $\tM$ is a well-ordering.
We obtain a representation in the form of either \eqref{PseudoReps} satisfying \eqref{PseudoRepCond}, or \eqref{FinalReduct} when $f=\el$.
However for the time being we just denote the multiplier in \eqref{PseudoReps} or \eqref{FinalReduct} as $\mar_\ipr$ since it satisfies $\mlt_\ipr (\mar_\ipr)=0$.

Finally we define the multiplier $\mar:=\gcd (\{\mar_\ipr\colon\ipr\in\cps\})$ such that $\gcd (\mar,\cpel)=1$ and $\mar=\sum_{\ipr\in\cps}\idr_\ipr\mar_\ipr$ with $\idr_\ipr\in\kux$.
By $\mar f=\sum_{\ipr\in\cps}\idr_\ipr\mar_\ipr f$ and the above representations of $\mar_\ipr f$ for $\ipr\in\cps$, we obtain the final conclusion.
\end{proof}

\begin{corollary}\label{Cor:LeadTermDef}
\begin{inparaenum}[(1)]
\item\label{item:ComptblDivide} The eliminant $\el$ is divisible by the compatible part $\cpel$ of the pre-eliminant $\pel$.

\item\label{item:ComptblProperBas} Let $I:=\langle\tas\rangle$ be the \zpi\ as in Lemma \ref{Lemma:PseudoReps} together with its pre-eliminant $\pel$ and prebasis $\tbs$.
When $\cpel\ne 1$, let us denote $\cpel$ simply as $\qr$ and define the PIR $\Rq$ that is isomorphic to $\kux/\pid\cpel$ as in \eqref{ProjectionPQR}.
Under the epimorphism $\mq$ as in \eqref{ExtendedProjectionPQR}, we denote $\Iq:=\mq (I)$ and $\Bq:=\mq (\tbs)$.
Then holds the identity \eqref{ProperLeadGen} that defines the \prb.
\end{inparaenum}
\end{corollary}
\begin{proof}
\begin{inparaenum}[(1)]
  \item A direct consequence of \eqref{FinalReduct}.

  \item Let $\ioq$ be the injection defined in \eqref{ExtendedEmbedding}.
  For every $\tf\in \Iq$ with $\lcc (\tf)\in\Rq^\ast$ being nonzero, $\exists\tg\in I$ such that $\mq (\tg)=\tf$.
  Then $\tg-\ioq (\tf)\in\langle\qr\rangle$ and thus $(\el/\qr)(\tg-\ioq (\tf))\in\langle\el\rangle$.
  Hence $(\el/\qr)\ioq (\tf)\in I$.

  As per Lemma \ref{Lemma:PseudoReps}, there exist $\{\qr_k\colon 0\le k\le s\}\subset\kxx$ as well as a multiplier $\mar\in (\kux)^\ast$ satisfying $\gcd (\mar,\qr)=1$ such that both \eqref{PseudoReps} and \eqref{PseudoRepCond} hold for $(\el/\qr)\ioq (\tf)$.
  For $1\le k\le s$, we collect the subscript $k$ into a set $\dss$ if it satisfies both $\lmc (\qr_kg_k)=\lmc (\ioq (\tf))$ and $\mq (\lcc (\qr_kg_k))\in\Rq^\ast$ being nonzero.
  Since $\mq (\mar\el/\qr)\in\Rq^\times$ is a unit, we have $\dss\ne\emptyset$ due to the equality \eqref{PseudoRepCond}.
  Then \eqref{ProperLeadGen} follows from the identity:
  $\ltc (\tf)=\mq (\mar\el/\qr)^{-1}\cdot\sum_{k\in\dss}\mq (\ltc (\qr_k))\cdot\ltc (\mq (g_k))\in\langle\ltc (\Bq)\rangle$.
\qedhere\end{inparaenum}
\end{proof}

When the compatible part $\cpel=1$, we simply ignore the partial proper basis $\Bq$ in the above Corollary \ref{Cor:LeadTermDef} \eqref{item:ComptblProperBas} and continue with the algorithm in the following Section \ref{Sec:ModularAlgm}.

\section{The Final Algorithmic Step in $\Rqx$}\label{Sec:ModularAlgm}

In this section we complete the construction of the \prb.
This is based on a modular algorithm which involves unorthodox computations in PIRs with zero divisors.
The content of this section totally parallels that of Section \ref{Sec:DivisionAlgm}.

Let $\pel$ be the pre-eliminant that is obtained in Lemma \ref{Algo:PseudoEliminant}.
Let $\cpel$ be its compatible part as before.
Suppose that $\pel/\cpel=\prod_{i=1}^s\tg_i^i$ is a squarefree factorization.
More specifically, $\tg_s\in\nonk$ is not a constant.
For $1\le i,j\le s$ and $i\ne j$, $\tg_i$ and $\tg_j$ are coprime if neither of them is a unit.
For $1\le i\le s$, if $\ipr\in\nonk$ is a monic irreducible polynomial that divides $\tg_i$, then $\mlt_\ipr (\tg_i)=1$ and $\mlt_\ipr (\pel/\cpel)=i$.

\begin{notation}[The modulus $\qr$ as the factor $\tg_i^i$ of $\pel/\cpel$]\label{Nott:ModulusQ}
\hfill

For $1\le i\le s$, if the squarefree factor $\tg_i$ of $\pel/\cpel$ is not a unit, let us simply denote the factor $\tg_i^i$ as $\qr$.
We define the PIR $\Rq$ that is isomorphic to $\kux/\pid{\tg_i^i}$ together with the epimorphism $\mq$ and injection $\ioq$ like in \eqref{ExtendedProjectionPQR} and \eqref{ExtendedEmbedding}.
\end{notation}

\begin{definition}[Modular proper term reduction in the algebra $\Rqx$]\label{Def:TermReductionPQR}
\hfill

For $f\in\Rqd$ and $g\in\rnd$ with $\lcc (g)\in\Rq^\ast$, suppose that $f$ has a term $\icr_\alpha\tx^\alpha$ with $\tx^\alpha\in\supp f\cap\langle\lmc (g)\rangle$.
We define the multipliers $\iur:=\mq (\lcm (l_\alpha,l_g)/l_\alpha)$ and $\lmr:=\mq (\lcm (l_\alpha,l_g)/l_g)$ with $l_\alpha:=\ioq (\icr_\alpha)$ and $l_g:=\ioq (\lcc (g))$.
If the multiplier $\iur$ satisfies $\iur\in\Rq^\times$ being a unit, we make a \emph{\mop\ reduction} of the term $\icr_\alpha\tx^\alpha$ by $g$ as follows.
\begin{equation}\label{TermReductionPQR}
h=\iur f-\frac{\lmr\tx^\alpha}{\lmc (g)}g.
\end{equation}
We call $h$ the \emph{remainder} of the reduction and $\iur$ the \emph{least multiplier} with respect to $g$.
\qed\end{definition}

\begin{definition}[Properly reduced polynomial in $\Rqx$]\label{Def:ProperlyReduced}
\hfill

A nonzero term $\icr_\alpha\tx^\alpha\in\Rqx$ is said to be \emph{properly reducible} with respect to $\tas=\{\Lst f1s\}\subset\Rqd$ if $\tx^\alpha$ is divisible by $\lmc (f_j)$ for some $f_j\in\tas$ and moreover, the least multiplier $\iur$ with respect to $f_j$ is a unit in $\Rq^\times$ as in \eqref{TermReductionPQR}.
A polynomial $f\in\Rqx$ is \emph{properly reduced} if it has no properly reducible terms with respect to $\tas$.
\qed\end{definition}

The proof for the following conclusion is almost a verbatim repetition of that for Theorem \ref{Thm:PseudoReduction}.

\begin{theorem}[Modular proper division or reduction in \text{$\Rqx$}]\label{Thm:ProperReduction}
\hfill

Suppose that $\tas=\{\Lst f1s\}$ is a finite polynomial set in $\Rqd$.
For every $f\in\Rqx$, there exist a multiplier $\mar\in\Rq^\times$ being a unit as well as a remainder $\rr\in\Rqx$ and quotients $q_j\in\Rqx$ for $1\le j\le s$ such that:
\begin{equation}\label{ProperDivisionExpression}
\mar f=\sum_{j=1}^s q_jf_j+\rr,
\end{equation}
where $\rr$ is properly reduced with respect to $\tas$, and the multiplier $\mar$ is literally a product of the least multipliers in \eqref{TermReductionPQR}.
Moreover, the polynomials in \eqref{ProperDivisionExpression} satisfy the following condition:
\begin{equation}\label{ProperDivisionCond}
\lmc (f)=\max\{\max_{1\le j\le s}\{\lmc (q_j)\cdot\lmc (f_j)\},\lmc (r)\}.
\end{equation}
\end{theorem}

For each pair $a,b\in\Rq$, let us define:
\begin{equation*}
\gcd\nolimits_\qr (a,b):=\mq (\gcd (\ioq (a),\ioq (b)));\quad
\lcm\nolimits_\qr (a,b):=\mq (\lcm (\ioq (a),\ioq (b))).
\end{equation*}

\begin{definition}[\Sp\ in $\Rqx$]\label{Def:SpolynPQR}
\hfill

Suppose that $f$ and $g$ are in $\Rqd$.
Let us denote $l_f:=\ioq(\lcc (f))$ and $l_g:=\ioq(\lcc (g))$.
We define the multipliers $\lmr_f:=\mq (\lcm (l_f,l_g)/l_f)$ and $\lmr_g:=\mq (\lcm (l_f,l_g)/l_g)$ as well as the monomial $\clm:=\lcm (\lmc (f),\lmc (g))$.
Then the following polynomial:
\begin{equation}\label{SPolyPQR}
S(f,g):=\frac{\lmr_f\clm}{\lmc (f)}f-\frac{\lmr_g\clm}{\lmc (g)}g
\end{equation}
is called the \emph{\Sp} of $f$ and $g$ in $\Rqx$.

There are two special \Sp s as follows.

\begin{inparaenum}[(1)]
\item\label{item:SpecSPolyQ} For every $f\in\Rqd$, we have $S(f,\qr)=\lnr_ff_1$ with $\lnr_f:=\mq (\qr/\gcd (l_f,\qr))$ and $f_1:=f-\ltc (f)$.

\item\label{item:SpecSPolyRq} For $f\in\Rqd$ and $\tg\in\rat$, we take $\lmc (\tg)=1$ and $\lcc (g)=g$ in \eqref{SPolyPQR}.
The \Sp\ satisfies $\idr S(f,\tg):=f_1g$ with $\idr=\gcd_\qr (\lcc (f),g)$ and $f_1:=f-\ltc (f)$.
\end{inparaenum}\qed
\end{definition}

The following two conclusions are modular versions of Lemma \ref{Lemma:RelativePrimePairs} and Lemma \ref{Lemma:TriangleIdentity} and we omit their trivial proofs here.

\begin{lemma}\label{Lemma:RelativePrimePQR}
For $f,g\in\Rqd$, suppose that $\lmc (f)$ and $\lmc (g)$ are relatively prime.
Let us denote $\idr:=\gcd_\qr (\lcc (f),\lcc (g))$.
Then the \Sp\ $S(f,g)$ satisfies $\idr S(f,g)=f_1g-g_1f$ with $f_1:=f-\ltc (f)$ and $g_1:=g-\ltc (g)$.
Moreover, we have $\lmc (S(f,g))=\max\{\lmc (f_1)\cdot\lmc (g),\lmc (g_1)\cdot\lmc (f)\}$.
\end{lemma}

For $f,g\in\rnd$ without both of them being in $\rat$, we define the \emph{relative common multiplier of $g$ versus $f$} as $\rcm (g\vert f):=\lmr_f\clm/\lmc (f)$ with $\lmr_f$ and $\clm$ being defined as in Definition \ref{Def:SpolynPQR}.

\begin{lemma}\label{Lemma:TriangleIdentityPQR}
For $f,g,h\in\rnd$ with at most one of them being in $\rat$, if $\lcm (\lmc (f),\lmc (g))$ is divisible by $\lmc (h)$, then holds the following triangular identity among the \Sp s:
\begin{equation}\label{TriangleIdentityPQR}
\mar S(f,g)=\frac{\mar\cdot\rcm (g\vert f)}{\rcm (h\vert f)}S(f,h)-\frac{\mar\cdot\rcm (f\vert g)}{\rcm (h\vert g)}S(g,h).
\end{equation}
Here the multiplier $\mar:=\mq (l_h/\idr)$ is nonzero, i.e., $\mar\in\Rq^\ast$ with $l_h:=\ioq (\lcc (h))$ and $\idr:=\gcd (\lcm (l_f,l_g),l_h)$.
\end{lemma}

The following algorithm is a slight revision of the one in Lemma \ref{Algo:PseudoEliminant}.

\begin{lemma}[The algorithm computing a partial \prb\ $\Bq$ in $\Rqx$]\label{Algo:ProperEliminant}
\hfill

Let $\tas\subset\Rqd$ be a finite polynomial set.
Then the algorithm that consists of the following sequence of instructions terminates in finite steps.
\begin{enumerate}
\renewcommand{\labelenumi}{\textbf{\upshape{B\theenumi}}}

  \item Start with a prebasis set $\tbs:=\tas$, temporary set $\CS:=\emptyset$ of \Sp s, and temporary pre-eliminant $\mel\in\Rq$ as $\mel:=0$.

  \item\label{item:AlgoStart} For each pair $\tf,\tg\in\tbs$ with $\tf\ne\tg$, execute the instructions in \textbf{\upshape{B\ref{item:SPolyCompMod}}}.

  \item\label{item:SPolyCompMod} If $\lmc (\tf)$ and $\lmc (\tg)$ are relatively prime, and the multiplier $\idr:=\gcd_\qr (\lcc (f),\lcc (g))$ is a unit, continue;

   Else if there exists an $h\in\tbs\setminus\{\tf,\tg\}$ such that $\lcm (\lmc (f),\lmc (g))\in\langle\lmc (h)\rangle$, and the multiplier $\mar$ as in \eqref{TriangleIdentityPQR} is a unit, continue;

   Else, compute the \Sp\ $S(\tf,\tg)$ as in \eqref{SPolyPQR} and add it into the set $\CS$.

   \item\label{item:SPolyReductPQR} For every \Sp\ $S\in\CS$, make a \mop\ reduction of $S$ by the prebasis set $\tbs$ to obtain a remainder that is denoted as $r$ like in Theorem \ref{Thm:ProperReduction}.
   Delete $S$ from the set $\CS$.

   If $\rr\in\Rqd$, add $\rr$ into $\tbs$.
   For every $\tf\in\tbs\setminus\{\rr\}$, execute the instructions in \textbf{\upshape{B\ref{item:SPolyCompMod}}} to
   compute the \Sp\ $S(\tf,\rr)$;

   Else if $\rr\in\rat$, redefine $\mel:=\mq (\gcd (\ioq (\rr),\qr))$ or $\mel:=\gcd_\qr (\rr,\mel)$, depending on $\mel=0$ or not.

   Else if $\rr\in\Rq^\times$ is a unit, halt the algorithm and output $\tbs=\{1\}$.

   \item\label{item:ReducSPoly} Repeat the instructions in \textbf{\upshape{B\ref{item:SPolyReductPQR}}} until $\CS=\emptyset$.

   \item\label{item:SpecialSPoly}
    For every $\tf\in\tbs$, depending on $\mel=0$ or not, compute the special \Sp\ $S(\tf,\qr)$ or $S(\tf,\mel)$ and add it into the set $\CS$ if $\mq(\gcd (l_f,\qr))$ or $\gcd_\qr (\lcc (\tf),\mel)$ is not a unit as in Definition \ref{Def:SpolynPQR} \eqref{item:SpecSPolyQ} or \eqref{item:SpecSPolyRq}.

    If the set $\CS$ is not empty, go back to \textbf{\upshape{B\ref{item:ReducSPoly}}}.

   \item Return the prebasis set $\tbs$, and the modulus $\qr$ or pre-eliminant $\mel$, depending on $\mel=0$ or not.

\end{enumerate}
\end{lemma}

\begin{lemma}\label{Lemma:nPQRSyzygy}
Let $\tas=\{f_j\colon 1\le j\le s\}\subset\Rqd$ be a finite polynomial set.
Suppose that for $1\le j\le s$, each $f_j$ has the same leading monomial $\lmc (f_j)=\tx^\alpha$.
\begin{enumerate}[(1)]
\item\label{item:SPolynomialExpansionPQR} If $f=\sum_{j=1}^s f_j$ satisfies $\lmc (f)\prec\tx^\alpha$, then
there exist nonzero multipliers $\ibr,\ibr_j\in\Rq^\ast$ for $1\le j<s$ such that
\begin{equation}\label{SPolynomialExpansionPQR}
\ibr f=\sum_{1\le j<s}\ibr_jS(f_j,f_s).
\end{equation}

\item\label{item:NondivisibilityPQR} For every irreducible factor $\ipr$ of the modulus $\qr$, we can always relabel the subscripts of the polynomial set $\tas=\{f_j\colon 1\le j\le s\}$ such that the multiplier $\ibr\in\Rq^\ast$ in \eqref{SPolynomialExpansionPQR} is not divisible by $\ipr$.
\end{enumerate}
\end{lemma}
\begin{proof}
\begin{inparaenum}[(1)]
\item Let us denote $l_j:=\ioq (\lcc (f_j))$ for $1\le j\le s$.
We define the multipliers $\lmr_j:=\lcm (l_j,l_s)/l_j$ for $1\le j<s$ and a multiplier $\aar:=\lcm_{1\le j<s}(\lmr_j)$. With $\aar_j:=\aar/\lmr_j\in\Rq^\ast$ for $1\le j<s$, we can easily verify that $\ibr_j:=\mq (\aar_j)$ and $\ibr:=\mq (\aar)$ in $\Rq^\ast$ satisfy the identity \eqref{SPolynomialExpansionPQR}.

\item Given an irreducible factor $\ipr$ of the modulus $\qr$, we can always change the order of the elements in the polynomial set $\tas=\{f_j\colon 1\le j\le s\}$ such that $\mlt_\ipr (l_s)=\min_{1\le j\le s}\{\mlt_\ipr (l_j)\}$.
Hence $\mlt_\ipr (\lmr_j)=\mlt_\ipr (l_s/\gcd (l_j,l_s))=0$ for $1\le j<s$, from which we can deduce that $\mlt_\ipr (\ibr)=0$.
\end{inparaenum}
\end{proof}

\begin{lemma}\label{Lemma:ModReps}
For a finite polynomial set $\tas\subset\Rqd$, let $\tbs=\{\tg_k\colon 1\le k\le s\}$ and $\mel$ be the prebasis and pre-eliminant that are obtained in the algorithm of Lemma \ref{Algo:ProperEliminant}.
For every $f\in\langle\tas\rangle$, there exist $\{\qr_k\colon 0\le k\le s\}\subset\Rqx$ and a unit multiplier $\mar\in \Rq^\times$ such that:
\begin{equation}\label{ProperReps}
\mar f=\sum_{k=1}^s\qr_k\tg_k+\qr_0\mel.
\end{equation}
Moreover, the polynomials in \eqref{ProperReps} satisfy the following condition:
\begin{equation}\label{ProperRepCond}
\lmc (f)=\max\bigl\{\max_{1\le k\le s}\{\lmc (\qr_k)\cdot\lmc (\tg_k)\},\lmc (\qr_0\mel)\bigr\}.
\end{equation}
In particular, the identities \eqref{ProperReps} and \eqref{ProperRepCond} still hold in the case of $\mel=0$.

Suppose that $I$ is a \zpi\ in $\kxx$ satisfying $\mq (I)=\langle\tas\rangle$.
Then its eliminant $\el$ also satisfies \eqref{ProperReps} with $\qr_0\in\Rq$ and $\mar\in \Rq^\times$ as follows:
\begin{equation}\label{FinalReductPQR}
\mar\mq (\el)=\qr_0\mel.
\end{equation}
\end{lemma}
\begin{proof}
The proof is almost a verbatim repetition of that for Lemma \ref{Lemma:PseudoReps}.
\end{proof}

\begin{corollary}
\begin{inparaenum}[(1)]
\item\label{item:UpdatedEliminant} In the case of $\mel=0$, the eliminant $\el$ is divisible by $\qr$;
Otherwise it is divisible by $\ioq (\mel)$.

\item Let $\tas\subset\Rqd$ be a finite polynomial set as in Lemma \ref{Lemma:ModReps} together with its pre-eliminant $\mel$ and prebasis $\tbs$.
Suppose that $I$ is a \zpi\ such that $\mq (I)=\langle\tas\rangle$.
In the case of $\mel=0$, we denote $\Iq:=\mq (I)$ and $\Bq:=\tbs$.
Then holds the identity \eqref{ProperLeadGen} that defines the \prb;
Otherwise we substitute $\ioq (\mel)$ for $\tg_i^i$ as the modulus $\qr$ in Notation \ref{Nott:ModulusQ}.
Then still holds the identity \eqref{ProperLeadGen} that defines the \prb.
\end{inparaenum}
\end{corollary}
\begin{proof}
\begin{inparaenum}[(1)]
\item The conclusion readily follows from the identity \eqref{FinalReductPQR}.

\item In the case of $\mel=0$, for every $\tf\in\Iq$, the identity \eqref{ProperReps} and \eqref{ProperRepCond} holds.
Then the identity \eqref{ProperLeadGen} follows from the nonempty subscript set $\dss:=\{1\le k\le s\colon\lmc (\qr_k)\cdot\lmc (\tg_k)=\lmc (\tf),\lcc (\qr_k)\cdot\lcc (\tg_k)\in\Rq^\ast\}$.

In the case when $\mel\in\Rq^\ast$ is nonzero, consider the natural epimorphism $\sigma\colon\Rq\rightarrow\Rq/\pid\mel$ with $\pid\mel$ being the principal ideal generated by $\mel\in\Rq$.
It induces an epimorphism $\sigma\colon\Rqx\rightarrow (\Rq/\pid\mel)[\tx]$ like in \eqref{ExtendedProjectionPQR}.
After we apply $\sigma$ to \eqref{ProperReps}, the proof of \eqref{ProperLeadGen} is a verbatim repetition of the above case.
\end{inparaenum}
\end{proof}

\section{An Example}
\label{Sec:Examples}

In the following example it is conspicuous that the coefficients of the \prb\ are of moderate sizes and swell like neither those of \gb\ basis over $\BQ$ nor those of \Mo's \gb\ basis over $\BQ [z]$.

\begin{example}
Suppose that we have a \zpi\ $I=\langle f,g,h\rangle\subset\BQ [x,y,z]$ with $f=-z^2(z+1)^3x+y;~g=z^4(z+1)^6x-y^2;~h=-x^2y+y^3+z^4(z-1)^5$.

For the purpose of comparison, let us list its reduced \gb\ basis $\tbs=\{\el,g_1,g_2,g_3,g_4\}$ over $\BQ$ with respect to the \lex\ ordering $x\succ y\succ z$ as follows:
\begin{align}
\el &=z^6(z-1)^5(z^{13}+9z^{12}+36z^{11}+84z^{10}+126z^9+126z^8+85z^7+31z^6+19z^5-9z^4+4z^3\notag\\
&\quad -4z^2-3z-1);\label{GrobnerEliminant}\\
g_1&=20253807z^2y+264174124z^{23}+1185923612z^{22}+850814520z^{21}-3776379304z^{20}\notag\\
&\quad -6824277548z^{19}+1862876196 z^{18}+12815317453z^{17}+3550475421z^{16}+2124010584z^{15}\notag\\
&\quad -35582561480z^{14}+42918431554z^{13}-41728834070z^{12}+35649844325z^{11}-17049238505z^{10}\notag\\
&\quad +3388659963z^9+930240431z^8-61146095z^7-518331181z^6;\notag\\
g_2&=20253807y^2+903303104z^{23}+4102316224z^{22}+3140448384z^{21}-12683487983z^{20}\notag\\
&\quad -23996669428z^{19}+4804720290z^{18}+43739947868z^{17}+14906482335z^{16}+9051639768z^{15}\notag\\
&\quad -121400613331z^{14}+139970660534z^{13}-138071007235z^{12}+118589702914z^{11}\notag\\
&\quad -55199680030z^{10}+11927452134z^9+2021069107z^8-38017822z^7-1768266833z^6;\notag\\
g_3&=2592487296z^2x+(7777461888z-2592487296)y+108083949263z^{23}+486376518055z^{22}\notag\\
&\quad +349557551130z^{21}-1558206505718z^{20}-2820179010211z^{19}+788268739077z^{18}\notag\\
&\quad +5350420983851z^{17}+1476923019345z^{16}+689330555757z^{15}-14602936038043z^{14}\notag\\
&\quad +17386123487861z^{13}-16350039201517z^{12}+13787524468420z^{11}-6235683207154z^{10}\notag\\
&\quad +786997920594z^9+628350552934z^8-64382649769z^7-206531133875z^6;\notag\\
g_4&=20253807x^2y+1037047036z^{23}+4686773132z^{22}+3455561112z^{21}-14868243976z^{20}\notag\\
&\quad -27470438972z^{19}+6731446644z^{18}+51651585868z^{17}+16267315284z^{16}+7429467573z^{15}\notag\\
&\quad -141636109619z^{14}+163168836472z^{13}-155454190640z^{12}+135706468958z^{11}\notag\\
&\quad -62903516282z^{10}+11263865469z^9+2500312823z^8+197272975z^7-1682438629z^6\allowdisplaybreaks[4]\notag\\
&\quad -101269035z^5+20253807z^4.\notag
\end{align}

According to \cite[P254, Theorem 4.5.12]{AL94}, this is also \Mo's strong \gb\ basis in $\BQ[z][x,y]$.

We begin the computation of the \prb\ by listing the prebasis elements $\tbs:=\{f,g,h\}$ in increasing order of their leading terms.

We disregard the \Sp\ $S(g,h)$ based on Lemma \ref{Lemma:TriangleIdentity} since $\lcm (\ltc (g),\ltc (h))=-z^4(z+1)^6x^2y$ is divisible by $\ltc (f)=-z^2(z+1)^3x$ and the multiplier $\mar=1$.
Then we compute the \Sp s $S(f,g)$ and $S(f,h)$ and make proper reductions.
The remainders are $e:=-y^2+z^2(z+1)^3y$ and $d:=z^2(z+1)^3[(z^4(z+1)^6-1)y+z^4(z-1)^5]$.
The prebasis set $\tbs$ becomes $\{d,e,f,g,h\}$.

We disregard the \Sp s $S(e,f)$ and $S(e,g)$ based on Lemma \ref{Lemma:RelativePrimePairs} since $\ltc (e)$ is relatively prime to both $\ltc (f)$ and $\ltc (g)$.
We also disregard the \Sp\ $S(e,h)$ since $\lcm (\lmc (e),\lmc (h))=x^2y^2$ is divisible by $\lmc (f)=x$.
Now the multiplier $\mar$ in the identity \eqref{TriangleIdentity} equals $-z^2(z+1)^3$ and we add it into the multiplier set $\mas$.

We disregard the \Sp s $S(d,f)$ and $S(d,g)$ since the leading monomial $\lmc (d)=y$ is relatively prime to $\lmc (f)=\lmc (g)=x$.
Moreover, we also disregard the \Sp\ $S(d,h)$ since $\lcm (\ltc (d),\ltc (h))=-z^2(z+1)^3(z^4(z+1)^6-1)x^2y$ is divisible by $\ltc (f)=-z^2(z+1)^3x$.

Let us compute and then make a proper reduction of the \Sp\ $S(d,e)=-z^4(z+1)^3(z^{13}+9z^{12}+36z^{11}+84z^{10}+126z^9+126z^8+85z^7+31z^6+19z^5-9z^4+4z^3
-4z^2-3z-1)y$.
We obtain the pre-eliminant $\pel=(z-1)^5z^8(z+1)^3(z^{13}+9z^{12}+36z^{11}+84z^{10}+126z^9+126z^8+85z^7+31z^6+19z^5-9z^4
+4z^3-4z^2-3z-1)$ with the multiplier $\iur=z^4(z+1)^6-1$.
We add $\iur$ into the multiplier set $\mas=\{z^2(z+1)^3,~z^4(z+1)^6-1\}$.

By a comparison between $\pel$ and $\mas$, we obtain the compatible part $\cpel=(z-1)^5(z^{13}+9z^{12}+36z^{11}+84z^{10}+126z^9+126z^8+85z^7+31z^6+19z^5-9z^4+4z^3
-4z^2-3z-1)$.

The squarefree factors of $\pel/\cpel$ are $\{z^8,(z+1)^3\}$.
For the modulus $\qr=(z+1)^3$ and hence over $\Rq=\BQ [z]/\pid{(z+1)^3}$, our computation shows that the ideal $\mq (I)=\{1\}$ and thus it should be disregarded.

For the modulus $\qr=z^8$ and hence over $\Rq=K[z]/\pid{z^8}$, we abuse the notations a bit and still use $f,g,h$ to denote $\mq (f),\mq (g),\mq (h)$.
The squarefree part of $\lcc (g)$ equals $z(z+1)$.
Hence we use it as a simple representation of $g$ and the same for $h$.
We disregard the \Sp\ $S(g,h)$ by the triangular identity with respect to $f$.

We compute and then make a \mop\ reduction of the \Sp s $S(f,g)$ and $S(f,h)$ to obtain the remainders $e:=-y^2+z^2(z+1)^3y$ and $d:=-z^2[(z+1)^3(z^4(z+1)^6+1)y-2z^5+z^4]$.
We add $d,e$ into the prebasis $\tbs=\{d,e,f,g,h\}$.

We disregard the \Sp s $S(e,f)$ and $S(e,g)$ since $\ltc (e)$ is relatively prime to both $\ltc (f)$ and $\ltc (g)$.
We also disregard the \Sp s $S(d,g)$ and $S(d,h)$ by the triangular identities with respect to $f$.
Moreover, the remainder of \mop\ reduction of the \Sp\ $S(d,e)$ equals $0$.

The \mop\ reduction of the \Sp\ $S(d,f)$ yields the pre-eliminant $\mel=z^6$.
We substitute $\qr=z^6$ for the old modulus $\qr=z^8$.
Then we compute and make a \mop\ reduction of the \Sp\ $S(e,h)$ to obtain $0$.
The same for the special \Sp s $S(f,\qr)$, $S(g,\qr)$ and $S(d,\qr)$.

After removing the redundancies of the prebasis $\tbs$, we obtain a partial \prb\ over $\Rq\cong\kux/\pid{z^6}$ as follows.
It corresponds to the pre-eliminant $\mel=z^6$ as above.
\begin{equation}\label{FirstPartialPropBas}
\Bq:~\biggl\{
\begin{aligned}
b_1&:=z^2(z+1)^3y;\\
b_3&:=z^2(z+1)^3x-y;
\end{aligned}\quad
\begin{aligned}
b_2&:=y^2;\\
b_4&:=x^2y-z^4(z-1)^5.
\end{aligned}
\biggr.
\end{equation}

The partial \prb\ over $\Rq\cong\kux/\pid{\cpel}$ is as follows.
It corresponds to the compatible part $\cpel$ that is also a factor of the eliminant $\el$ by Corollary \ref{Cor:LeadTermDef} \eqref{item:ComptblDivide}.
\begin{equation}\label{SecondPartialPropBas}
\Bq:~\biggl\{
\begin{aligned}
a_1&:=z^2(z+1)^3(z^4(z+1)^6-1)y+z^6(z+1)^3(z-1)^5;\\
a_2&:=z^4(z+1)^6(z^4(z+1)^6-1)x+z^6(z+1)^3(z-1)^5.
\end{aligned}
\biggr.
\end{equation}

Hence the eliminant is $\el=\cpel\cdot\ioq (\mel)$, which coincides with the one obtained via the \gb\ basis in \eqref{GrobnerEliminant}.
The partial proper bases in \eqref{FirstPartialPropBas} and \eqref{SecondPartialPropBas} constitute the \prb\ of $I$ as in \eqref{NewBases}.
\end{example}

\section{Benchmark Testings in a Cascade of Complexity}
\label{Sec:Timings}

In order to corroborate the distinctive computational advantage of the \prb\ algorithm over \bg's and \Mo's classical algorithms for \gb\ bases, we conducted a cascade of benchmark testings that are in increasing order of complexity.
The programming language is Maple 18.
The timings were conducted on an 11th Gen Intel Core i7 3.30 GHz system with 32.0 GB RAM under a 64-bit version of Windows 10 operating system.
All the computations are over the rationals $\BQ$ with respect to the typical \lex\ ordering $x\succ y\succ z$.

\begin{example}\label{Expl:Kitties}
\begin{enumerate}[(1)]
\item $I=\langle x^2+x-zy,~-zx+y^3+2,~-x+y+z^2-1\rangle$;
\item $I=\langle x^2+x-zy^2,~-zx+y^4+2,~-x+y^2+z^2-1\rangle$;
\item $I=\langle x^2+x-zy^3,~-zx+y^4+2,~-x+y^2+z^2-1\rangle$;
\item $I=\langle x^2+x-zy^3,~-zx+y^5+2,~-x+y^2+z^2-1\rangle$;
\item $I=\langle x^2+x-zy^3,~-z^2x+y^5+2,~-x+y^2+z^3-1\rangle$;
\item $I=\langle x^2+x-z^2y^3,~-z^3x+y^5+2,~-x+y^2+z^3-1\rangle$;
\item $I=\langle x^2+x-zy^3,~-zx+y^5+2,~-x+y^3+z^2-1\rangle$;
\item $I=\langle x^3+x-zy^3,~-zx+y^5+2,~-x+y^3+z^2-1\rangle$;
\item $I=\langle x^3+x-z^2y^4,~-z^2(z^3-2)x^2+y^6+2z^4,~-x^2+y^4+z^4-z\rangle$;
\item $I=\langle x^2+x-z^2y^3,~-z(z-1)x^2+y^6+z^4,~-x+y^5+z^2-1\rangle$;
\item $I=\langle x^3+x-z^2y^4,~-z(z-1)x^2+y^5+z^4,~-x+y^5+z^2-1\rangle$;
\item $I=\langle x^2+x-z^2y^4,~-z(z-1)x^2+y^6+z^4,~-x+y^5+z^2-1\rangle$;
\item $I=\langle x^3+x-z^2y^4,~-z(z-1)x^2+y^6+z^4,~-x+y^5+z^2-1\rangle$;
\end{enumerate}
\end{example}

In the following tables the unit for the timings is in second.
The symbol ``N/A" means that either Maple crashed or the running time was more than $86400.000$ seconds (24 hours).

\begin{longtable}{|c|r|r|r|}
\hline
Example \ref{Expl:Kitties} & \multicolumn{1}{|c|}{Proper basis} & \multicolumn{1}{|c|}{\Mo's basis} & \multicolumn{1}{|c|}{\bg's basis}\\\hline\allowdisplaybreaks[4]
(1)  &  $0.016$  &  $0.032$  & $0.453$ \\\hline
(2)  &  $0.000$  &  $0.016$  & $0.047$ \\\hline
(3)  &  $0.750$   &  $1.594$  & $79.250$ \\\hline
(4)  &  $2.172$   &  $7.547$  & $442.125$ \\\hline
(5)  &  $1.953$   &  $112.391$  & $9569.672$ \\\hline
(6)  &  $2.000$   &  $145.703$  & $14694.797$ \\\hline
(7)  &  $2.000$   &  $36229.156$  & N/A \\\hline
(8)  &  $2.953$   &  $12516.609$  & N/A \\\hline
(9)  &  $5.829$   &  N/A & N/A \\\hline
(10)  &  $15.390$   &  N/A & N/A \\\hline
(11)  &  $32.203$   &  N/A & N/A \\\hline
(12)  &  $75.297$   &  N/A & N/A \\\hline
(13)  &  $194.735$   &  N/A & N/A \\\hline
\end{longtable}

\begin{remark}
The above implementations for all the three algorithms are primitive in the sense that we only make basic optimizations like in Lemma \ref{Lemma:RelativePrimePairs} and Lemma \ref{Lemma:TriangleIdentity}, or Lemma \ref{Lemma:RelativePrimePQR} and Lemma \ref{Lemma:TriangleIdentityPQR}.
This leaves room for further improvements and optimizations on the \prb\ algorithm.
In fact, it is a natural expectation that all the state of the art algorithms that improves \bg's algorithm over the last decades can be applied to the proper basis for further improvements and optimizations.
\end{remark}

\section{Conclusion and Acknowledgement}
\label{Sec:Conclusion}

We defined the \prb\ based on the proper division in Definition \ref{Def:NewBases} for \zpi s.
The \prb\ algorithm is more efficient than \bg's and \Mo's classical algorithms for \gb\ bases.
This is corroborated by the benchmark testings with respect to the typical \lex\ ordering over the rationals $\BQ$ in Section \ref{Sec:Timings}.
We shall address the significant progress in generalizing the \prb\ to the polynomial ideals of positive dimensions.
It shall be interesting to make further efficiency improvements on the \prb\ algorithm in the future.

%


\begin{thebibliography}{199}

\bibitem{AL94}
Adams, W., Loustaunau, P., 1994. An Introductin to \gb\ Bases. Grad. Stud. Math. 3, Amer. Math. Soc.

\bibitem{Arn03}
Arnold, E., 2003. Modular algorithms for computing \gb\ bases. J. Symbolic Comput. 35, 403-419.

\bibitem{BW93}
Becker, T., Weispfenning, V., 1993. \gb\ Bases. A computational approach to commutative algebra. Grad. Texts in Math. 141, Springer.

\bibitem{BW98}
\bg, B., Winkler, F., 1998. \gb\ Bases and Applications. London Math. Soc. Lecture Note Ser. 251, Cambridge Univ. Press.

\bibitem{Buc65}
\bg,\ B., 2006. 1965 Ph.D. Thesis: An algorithm for finding the basis elements of the residue class ring of a
zero dimensional polynomial ideal. J. Symbolic Comput. 41 (3-4), 475-511.

\bibitem{Buc85}
\bg,\ B., 1985. \gb\ bases: an algorithmic method in polynomial ideal theory, in: Bose, N. (Eds.), Multidimensional Systems Theory. D. Reidel Publishing, pp. 184-232.

\bibitem{CKM97}
Collart, S., Kalkbrener, M., Mall, D., 1997. Converting bases with the \gb\ walk. J. Symbolic Comput. 24, 465-469.

\bibitem{CLO05}
Cox, D., Little, J., O'Shea, D., 2005. Using Algebraic Geometry, second ed. Grad. Texts in Math. 185, Springer-Verlag.

\bibitem{CLO15}
Cox, D., Little, J., O'Shea, D., 2015. Ideals, Varieties, and Algorithms. An introduction to computational algebraic geometry and commutative algebra, fourth ed. Springer-Verlag.

\bibitem{DL06}
Decker, W., Lossen, C., 2006. Computing in Algebraic Geometry. Springer-Verlag.

\bibitem{Dub90}
Dub\'e, T., 1990. The structure of polynomial ideals and \gb\ bases. SIAM J. Comput, 19 (4), 750-773.

\bibitem{EF17}
Eder, C., Faug\`ere, J., 2017. A survey on signature-based algorithms for computing \gb\ bases. J. Symbolic Comput. 80, 719-784.

\bibitem{EH12}
Ene, V., Herzog, J., 2012. \gb\ Bases in Commutative Algebra. Grad. Stud. Math. 130, Amer. Math. Soc.

\bibitem{EH21}
Eder, C., Hofmann, T., 2021. Efficient \gb\ bases computation over principal ideal rings. J. Symbolic Comput. 103, 1-13.

\bibitem{Ebe83}
Ebert, G., 1983. Some comments on the modular approach to \gb-bases. ACM SIGSAM Bulletin 17, 28-32.

\bibitem{FGLM93}
Faug\`ere, J., Gianni, P., Lazard, D., Mora, T., 1993. Efficient computation of zero-dimensional \gb\ bases by change of ordering.
J. Symbolic Comput. 16, 329-344.

\bibitem{FV20}
Francis, M., Verron, T., 2020. A signature-based algorithm for computing \gb\ bases over \Pid s. Math. Comput. Sci. 14, 515-530.

\bibitem{Fau02}
Faug\`ere, J., 2002. A new efficient algorithm for computing \gb\ bases without reduction to zero (F5), in: ISSAC 2002, Proceedings of the 2002 International Symposium on Symbolic and Algebraic Computation, ACM Press, 75-83.

\bibitem{Fro97}
Fr\"oberg, R., 1997. An Introduction to \gb\ Bases. John Wiley \& Sons.

\bibitem{GG13}
von zur Gathen, J., Gerhard, J., 2013. Modern Computer Algebra, third ed. Cambridge Univ. Press.

\bibitem{GM88}
Gebauer, R., \Mo, M., 1988. On an installation of \bg's algorithm. J. Symbolic Comput. 6, 275-286.

\bibitem{GMN91}
Giovini, A., Mora, T., Niesi, G., Robbiano, L., Traverso, C., ``One sugar cube, please," or selection strategies in the buchberger algorithm, in: Watt, S. (Eds.), ISSAC 1991, Proceedings of the 1991 International Symposium on Symbolic and Algebraic Computation, ACM Press, 49-54.

\bibitem{GP08}
Greuel, G., Pfister, G., 2008. A Singular Introduction to Commutative Algebra. Springer-Verlag.

\bibitem{GVW16}
Gao, S.H., Volny IV, F., Wang, M.S., 2016. A new framework for computing \gb\ bases. Math. Comp. 85 (297), 449-465.

\bibitem{Gra93}
Gr\"abe, H., 1993. On lucky primes. J. Symbolic Comput. 15, 199-209.

\bibitem{HH22}
Harvey, D., van der Hoeven, J., 2022. Polynomial multiplication over finite fields in time $O(n\log n)$. J. ACM. 69 (2), 12:1-12:40.

\bibitem{KR00}
Kreuzer, M., Robbiano, L., 2000. Computational Commutative Algbera 1. Springer-Verlag.

\bibitem{Laz92}
Lazard, D., 1992. Solving zero-dimensional algebraic systems. J. Symbolic Comput. 13 (2), 117-131.

\bibitem{Lic12}
Lichtblau, D., 2012. Effective computation of strong \gb\ bases over euclidean domains. Illinois J. Math. 56 (1), 177-194.

\bibitem{Mis93}
Mishra, B., 1993. Algorithmic Algebra. Texts Monogr. Comput. Sci. Springer-Verlag.

\bibitem{Mol88}
\Mo, H., 1988. On the construction of \gb\ bases using syzygies. J. Symbolic Comput. 6 (2), 345-359.

\bibitem{NS01}
Norton, G., S\u{a}l\u{a}gean, A., 2001. Strong \gb\ bases for polynomials over a principal ideal ring. Bull. Aust. Math. Soc. 64, 505-528.

\bibitem{Nab09}
Nabeshima, K., 2009. Reduced \gb\ bases in polynomial rings over a polynomial ring. Math. Comput. Sci. 2, 587-599.

\bibitem{Pau92}
Pauer, F., 1992. On lucky ideals for \gb\ basis computations. J. Symbolic Comput. 14, 471-482.

\bibitem{Rou99}
Rouillier, F., 1999. Solving zero-dimensional systems through the rational univariate representation. Appl. Algebra Engrg. Comm. Comput. 9, 433-461.

\bibitem{SS06}
Suzuki, A., Sato, Y., 2006. A simple algorithm to compute comprehensive \gb\ bases using \gb\ bases, in: ISSAC 2006, Proceedings of the 2006 International Symposium on Symbolic and Algebraic Computation, ACM Press, 326-331.

\bibitem{ST89}
Sasaki, T., Takeshima, T., 1989. A modular method for \gb-basis construction over $\BQ$ and solving system of algebraic equations. J. Information Processing 12, 371-379.

\bibitem{Stu95}
Sturmfels, B., 1995. \gb\ Bases and Convex Polytopes. Univ. Lecture Ser. 8, Amer. Math. Soc.

\bibitem{Tra89}
Traverso, C., 1989. \gb\ trace algorithms, in: Symbolic and Algebraic Computations (Rome 1988). Lecture Notes in Comput. Sci. 358, Springer-Verlag, 125-138.

\bibitem{Win88}
Winkler, F., 1988. A $p$-adic approach to the computation of \gb\ bases. J. Symbolic Comput. 6, 287-304.

\bibitem{Wu83}
Wu, W.T., 1983. On the decision problem and the mechanization of theorem-proving in elementary geometry, in: Automated Theorem Proving: After 25 Years. Bledsoe, W., Loveland, D. (Eds.) Contemp. Math. 29, Amer. Math. Soc., 213-234.

\bibitem{Wu01}
Wu, W.T., 2000. Mathematics Mechanization: Mechanical Geometry Theorem-Proving, Mechanical Geometry Problem-Solving and Polynomial Equations-Solving. Math. Appl. 489, Kluwer Dordrecht and Science Press Beijing.

\end{thebibliography}
\end{document}